\documentclass[a4paper,oneside,reqno,11pt]{article}
\usepackage{latexsym,amsthm,amscd,amsmath,amssymb}
\usepackage[mathscr]{eucal}
\usepackage[english]{babel}

\newtheorem{e-proposition}[theorem]{Proposition}

\newtheorem{e-definition}[theorem]{Definition\rm}

\def\n{\noindent}
\def\intn{\int\nolimits}

\def\ul#1{{\smash{\rm#1}}}
\def\SUPP{\mathop{\ul{SUPP}}}
\def\IM{\mathop{\ul{IM}}}
\def\IC{\mathop{\ul{IC}}}
\def\IT{\mathop{\ul{IT}}}
\def\bu{{\,{\scriptscriptstyle\bullet}\,}}

\begin{document}

\centerline{\bf The weighted log canonical threshold}
\vskip4mm
\centerline{\bf Pham Hoang Hiep}
\vskip1cm

{\bf Abstract:} In this note, we show how to apply the original $L^2$-extension theorem of Ohsawa and Takegoshi to the standard basis of a multiplier ideal sheaf associated with a plurisubharmonic function. In this way, we are able to reprove the strong openness conjecture and to obtain an effective version of the semicontinuity theorem for weighted log canonical thresholds.

\section{Introduction and main results}
\vskip-2mm

Let $\Omega$ be a domain in $\mathbb C^n$ and $\varphi$ in the set $\hbox{\rm PSH}(\Omega)$ of plurisubharmonic functions on~$\Omega$. Following Demailly and Koll\'ar \cite{DK00}, we introduce the log canonical threshold of $\varphi$ at a point $z_0\in\Omega$
$$c_\varphi (z_0) = \sup\big\{ c>0:\ e^{-2c\,\varphi} \hbox{ is } L^1 \hbox{ on a neighborhood of }z_0\big\}\in (0,+\infty].$$
It is an invariant of the singularity of $\varphi$ at $z_0$. We refer to \cite{Ca12}, \cite{De93a}, \cite{De93b}, \cite{FEM03}, \cite{DH12}, \cite{DK00}, \cite{FEM10} \cite{HHH12}, \cite{Hi13}, \cite{Ki94}, \cite{Na90}, \cite{PS00}, \cite{Sk72} for further information about this number. In \cite{DK00}, Demailly and Koll\'ar stated the following openness conjecture.\vskip3pt
\n
{\bf Conjecture.} {\em The set}~~$\{c>0:\ e^{-2c\,\varphi}~\hbox{\em is }~L^1~\hbox{\em on a neighborhood of }z_0\big\}$ {\em equals the open interval}~~$(0,c_\varphi (z_0)).$
\vskip3pt
\n
In 2005, this conjecture was proved in dimension $2$ by Favre and Jonsson (\cite{FJ05}, \cite{JM12}, \cite{JM14}). In 2013, Berndtsson (\cite{Be13}) completely proved it in arbitrary dimension. For every holomorphic function $f$ on $\Omega$, we introduce the {\em weighted log canonical threshold} of $\varphi$ with weight $f$ at $z_0$:
$$c_{\varphi, f} (z_0) = \sup\big\{ c>0:\ |f|^2 e^{-2c\,\varphi} \hbox{ is } L^1 \hbox{ on a neighborhood of }z_0\big\}\in (0,+\infty].$$
Recently, Guan-Zhou used a sophisticated version of the $L^2$-extension theorem of Ohsawa and Takegoshi in combination with the curve selection lemma, to prove the ``strong'' openness conjecture, i.e.\ the analogue openness statement for weighted thresholds $c_{\varphi, f}(z_0)$, and a related semi-continuity theorem for the weighted log canonical threshold (\cite{GZ12}, \cite{GZ14}). In this note, we show how one can apply the original version \cite{OT87} of the $L^2$-extension theorem to the members of a standard basis for a multiplier ideal sheaf of holomorphic functions associated with a plurisubharmonic function $\varphi$. In this way, by means of a simple induction on dimension, we reprove the strong openness conjecture, and give an effective version of the semicontinuity theorem for weighted log canonical thresholds. The main results are contained in the following theorem.
\vskip3pt
\n
{\bf Main theorem.} {\it Let $f$ be a holomorphic function on an open set
$\Omega$ in ${\mathbb C}^n$ and let $\varphi\in\hbox{\rm PSH}(\Omega)$.}
\vskip2pt
\n
(i) (``Semicontinuity theorem'') {\it Assume that $\intn_{\Omega'}e^{-2c\,\varphi}dV_{2n}<+\infty$ on some open subset $\Omega'\subset\Omega$ and let $z_0\in\Omega'$. Then for
$\psi\in\hbox{\rm PSH}(\Omega')$, there exists $\delta=\delta(c,\varphi,\Omega',z_0)>0$ such that $\Vert \psi-\varphi\Vert_{L^1(\Omega')}\le \delta$ implies \hbox{$c_{\psi}(z_0)>c$}. Moreover, as $\psi$ converges to $\varphi$ in~$L^1(\Omega')$, the function $e^{-2c\,\psi}$ converges to $e^{-2c\,\varphi}$ in~$L^1$ on every relatively compact open subset $\Omega''\Subset\Omega'$.}
\vskip2pt
\n
(ii) (``Strong effective openness'') {\it Assume that $\intn_{\Omega'}|f|^2e^{-2c\,\varphi}dV_{2n}<+\infty$ on some open subset $\Omega'\subset\Omega$. When $\psi\in\hbox{\rm PSH}(\Omega')$ converges to $\varphi$ in $L^1(\Omega')$ with $\psi\le\varphi$, the function $|f|^2e^{-2c\,\psi}$ converges to $|f|^2e^{-2c\,\varphi}$ in $L^1$ norm on every relatively compact open subset $\Omega''\Subset\Omega'$.}
\vskip3pt
\n
{\bf Corollary 1.1} {\rm(``Strong openness'')}. {\it For any plurisubharmonic function $\varphi$ on a neighborhood of a point $z_0\in{\mathbb C}^n$, the set
$$\{ c>0:\ |f|^2 e^{-2c\,\varphi} \hbox{ is } L^1 \hbox{ on a neighborhood of }z_0\}$$
is an open interval $(0,c_{\varphi,f}(z_0))$.}
\vskip3pt
\n
{\bf Corollary 1.2} {\rm(``Convergence from below'')}. {\it If $\psi\le\varphi$ converges to $\varphi$ in a neighborhood of $z_0\in{\mathbb C}^n$, then $c_{\psi, f}(z_0)\leq c_{\varphi, f}(z_0)$ converges to $c_{\varphi, f}(z_0)$.}
\vskip3pt
\n
In fact, after subtracting a large constant to~$\varphi$, we can assume $\varphi\leq 0$ in both corollaries. Then Cor.~1.1 is a consequence of assertion (ii) of the main theorem when we take $\Omega'$ small enough and $\psi=(1+\delta)\varphi$ with~$\delta\searrow 0$. In Cor.~1.2, we have by definition $c_{\psi,f}(z_0)\leq c_{\varphi,f}(z_0)$ for $\psi\leq\varphi$, but again (ii) shows that $c_{\psi,f}(z_0)$ becomes${}\geq c$ for any given value $c\in(0,c_{\varphi,f}(z_0))$, whenever $\Vert\psi-\varphi\Vert_{L^1(\Omega')}$ is sufficiently small.
\vskip3pt
\n
{\bf Remark 1.3.} One cannot remove condition $\psi\leq\varphi$ in assertion (ii) of the main theorem. Indeed, let us choose $f(z)=z_1$, $\varphi(z)=\log |z_1|$ and $\varphi_j(z) = \log |z_1+\frac {z_2} j|$, for $j\geq 1$. We have $\varphi_j\to\varphi$ in $L_{\rm loc}^1({\mathbb C}^n)$, however $c_{\varphi_j, f}(0)=1<c_{\varphi,f}(0)=2$ for all $j\geq 1$. On the other hand, condition (i) does not require any given inequality between $\varphi$ and $\psi$. Modulo Berndtsson's solution of the openness conjecture, (i) follows from the effective semicontinuity result of \cite{DK00}, but (like Guan and Zhou) we reprove here both by a direct and much easier method.
\vskip3pt
\n
{\bf Remark 1.4.} As in Guan-Zhou \cite{GZ13a}, \cite{GZ13b}, one can reformulate Cor.~1.1 in terms of multiplier ideal sheaves. Denote by ${\mathcal I}(c\,\varphi)$ the sheaf of germs of holomorphic functions $f\in{\mathcal O}_{{\mathbb C}^n,z}$ such that $\int_U|f|^2e^{-2c\,\varphi}dV_{2n}<+\infty$ on some neighborhood $U$ of $z$ (it is known by \cite{Na90} that this is a coherent ideal sheaf over~$\Omega$, but we will not use this property here). Then at every point $z\in\Omega$ we have
$${\mathcal I}(c\,\varphi)={\mathcal I}_+(c\,\varphi):=\lim_{\epsilon\searrow 0}\,
{\mathcal I}((1+\epsilon)c\,\varphi).$$

\section{Proof of the main theorem}
\vskip-2mm
We equip the ring ${\mathcal O}_{\mathbb C^n,0}$ of germs of holomorphic functions at $0$ with the homogeneous lexicographic order of monomials $z^\alpha = z_1^{\alpha_1} \ldots z_n^{\alpha_n}$, that is, $z_1^{ \alpha_1 } \ldots z_n^{ \alpha_n }<z_1^{ \beta_1 } \ldots z_n^{ \beta_n }$ if and only if $| \alpha |=\alpha_1+\ldots+\alpha_n<| \beta |=\beta_1+\ldots+\beta_n$ or $| \alpha |=| \beta |$ and $\alpha_i < \beta_i$ for the first index $i$ with $\alpha_i \not = \beta_i$. For each $f(z)=a_{\alpha^1} z^{\alpha^1}+a_{\alpha^2} z^{\alpha^2}+\ldots\;$ with $a_{\alpha^j}\not = 0$, $j\geq 1$ and $z^{ \alpha^1 } < z^{ \alpha^2 } < \ldots\;$, we define the {\it initial coefficient}, {\it initial monomial} and  {\it initial term} of $f$ to be respectively $\IC(f)=a_{\alpha^1}$, $\IM(f)=z^{\alpha^1}$, $\IT(f)=a_{\alpha^1} z^{\alpha^1}$, and the support of $f$ to be $\SUPP(f)=\{z^{ \alpha^1 }, z^{ \alpha^2 }, \ldots\}$. For any ideal $\mathcal I$~of~${\mathcal O}_{\mathbb C^n,0}$, we define $\IM(\mathcal I)$ to be the ideal generated by $\{ \IM(f) \}_{\{ f\in\mathcal I\} }$. First, we recall the division theorem of Hironaka and the concept of standard basis of an ideal.
\vskip3pt
\n
{\bf Division theorem of Hironaka} (see \cite{Ga79}, \cite{Ba82}, \cite{BM87}, \cite{BM89}, \cite{Ei95}). {\it Let $f, g_1,\ldots,g_k\in {\mathcal O}_{\mathbb C^n,0}$. Then there
exist $h_1,\ldots,h_k,s\in {\mathcal O}_{\mathbb C^n,0}$ such that
$$f=h_1g_1+\ldots+h_kg_k+s,$$
and $\SUPP(s)\cap \langle\IM(g_1),\ldots,\IM(g_k)\rangle =\emptyset$, where $\langle\IM(g_1),\ldots,\IM(g_k)\rangle$ denotes the ideal generated by the family~$(\IM(g_1),\ldots,\IM(g_k))$.}
\vskip3pt
\n
{\bf Standard basis of an ideal.} Let $\mathcal I$ be an ideal of ${\mathcal O}_{\mathbb C^n,0}$ and let $g_1,\ldots,g_k\in\mathcal I$ be such that $\IM(\mathcal I)=\langle\IM(g_1),\ldots,\IM(g_k)\rangle$. Take $f\in\mathcal I$. By the division theorem of Hironaka, there exist $h_1,\ldots,h_k,s\in {\mathcal O}_{\mathbb C^n,0}$ such that
$$f=h_1g_1+\ldots+h_kg_k+s,$$
and $\SUPP(s)\cap\IM(\mathcal I) =\emptyset$. On the other hand, since $s=f-h_1g_1+\ldots+h_kg_k\in\mathcal I$, we have \hbox{$\IM(s)\in\IM(\mathcal I)$}. Therefore $s=0$ and the $g_j$'s are generators of~$\mathcal I$. By permuting the $g_j$'s and performing ad hoc subtractions, we can always arrange that
$\IM(g_1)<\IM(g_2)<\ldots<\IM(g_k)$, and we then say
that $(g_1,\ldots,g_k)$ is a standard basis of $\mathcal I$.
\vskip5pt
\n
We will prove the main theorem by induction on dimension $n$. Of course, it holds for $n=0$. Assume that the theorem holds for dimension $n-1$. Thanks to the $L^2$-extension theorem of Ohsawa and Takegoshi (\cite{OT87}), we obtain the following key lemma.
\vskip3pt
\n
{\bf Lemma 2.1.} {\it Let $\varphi\leq 0$ be a plurisubharmonic function and $f$ be a holomorphic function on the polydisc $\Delta_R^n$ of center $0$ and $($poly$)$radius $R>0$ in $\mathbb C^n$, such that for some $c>0$
$$\intn_{\Delta_R^n} |f(z)|^2 e^{-2c\,\varphi(z)} d V_{2n}(z)<+\infty.$$
Let $\psi_j\leq 0$, $j\geq 1$, be a sequence of plurisubharmonic functions on
$\Delta_R^n$ with $\psi_j\to\varphi$ in $L^1_{\rm loc}(\Delta_R^n)$, and assume
that either $f=1$ identically or $\psi_j\leq\varphi$ for all~$j\geq 1$.
Then for every $r<R$ and $\epsilon\in(0,\frac{1}{2}r]$, there exist a value $w_n\in\Delta_\epsilon \smallsetminus\{0\}$, an index $j_0$, a constant $\tilde c>c$  and a sequence of holomorphic functions $F_j$ on $\Delta_{r}^n$, $j\geq j_0$, such that $\IM(F_j)\leq \IM(f)$, $F_j(z) = f(z) + (z_n-w_n)\sum a_{j,\alpha} z^{\alpha}$ with $| w_n| | a_{j,\alpha} |\leq r^{-|\alpha|}\epsilon$ for~all~$\alpha\in{\mathbb N}^n$, and
$$\intn_{ \Delta_{r}^n } |F_j(z)|^2 e^{-2\tilde c\,\psi_j(z)} d V_{2n}(z)\leq\frac{\epsilon^2}{|w_n|^2}<+\infty,~~~\forall j\geq j_0.$$
Moreover, one can choose $w_n$ in a set of positive measure in the
punctured disc $\Delta_\epsilon\smallsetminus\{0\}$ $($the index $j_0=j_0(w_n)$ and the constant $\tilde c=\tilde c(w_n)$ may then possibly depend on $w_n)$.}
\vskip3pt
\n
{\it Proof.} By Fubini's theorem we have
$$\intn_{\Delta_R}\bigg[\intn_{\Delta_R^{n-1}} |f(z',z_n)|^2 e^{-2c\,\varphi(z',z_n)} d V_{2n-2}(z')\bigg] dV_2(z_n)<+\infty.$$
Since the integral extended to a small disc $z_n\in \Delta_\eta$ tends to $0$ as $\eta\to 0$, it will become smaller than any preassigned value, say $\epsilon_0^2>0$, for $\eta\leq\eta_0$ small enough. Therefore we can choose a set of positive measure of values $w_n\in\Delta_\eta\smallsetminus \{ 0\}$ such that
$$\intn_{\Delta_R^{n-1}} |f(z',w_n)|^2e^{-2c\,\varphi(z',w_n)} d V_{2n-2}(z')\leq\frac { \epsilon_0^2}{\pi\eta^2}<\frac {\epsilon_0^2}{|w_n|^2}.$$
Since the main theorem is assumed to hold for $n-1$, for any $\rho<R$ there exist $j_0=j_0(w_n)$ and $\tilde c=\tilde c(w_n)>c$ such that
$$\intn_{\Delta_{\rho}^{n-1}} |f (z',w_{n}) |^2 e^{-2\tilde c\,\psi_j(z',w_n)} d V_{2n-2}(z')<\frac {\epsilon_0^2}{|w_n|^2},~~~\forall j\geq j_0.$$
(For this, one applies part (i) in case $f=1$, and part (ii) in case $\psi_j\leq\varphi$, using the fact that $\psi=\frac{\tilde c}{c}\,\psi_j$ converges to $\varphi$ as $\tilde c\to c$ and $j\to+\infty$). Now, by the $L^2$-extension theorem of Ohsawa and Takegoshi (see \cite{OT87} or \cite{De09}), there exists a holomorphic function $F_j$ on $\Delta_{\rho}^{n-1}\times\Delta_R$ such that $F_j(z',w_n)=f(z',w_n)$ for all $z'\in\Delta_{\rho}^{n-1}$,
and
\begin{eqnarray*}
&&\intn_{\Delta_{\rho}^{n-1}\times \Delta_R} |F_j(z)|^2 e^{-2\tilde c\,\psi_j(z)} d V_{2n}(z)\\
&\leq& C_nR^2\intn_{\Delta_{\rho}^{n-1}} |f(z',w_n)|^2 e^{-2\tilde c\,\psi_j(z',w_n)} d V_{2n-2}(z')\\
&\leq&\frac { C_n R^2 \epsilon_0^2}{ |w_n|^2 },
\end{eqnarray*}
where $C_n$ is a constant which only depends on $n$ (the constant is universal for $R=1$ and is rescaled by $R^2$ otherwise). By the mean value inequality for the plurisubharmonic function $|F_j|^2$, we get
\begin{eqnarray*}
|F_j(z)|^2&\leq& \frac 1 {\pi^n(\rho -|z_1|)^2\ldots(\rho -|z_n|)^2}\intn_{ \Delta_{\rho -|z_1|} (z_1)\times\ldots\times\Delta_{\rho -|z_n|} (z_n) } |F_j|^2 dV_{2n}\\
&\leq& \frac {C_n R^2 \epsilon_0^2}{\pi^n(\rho -|z_1|)^2\ldots(\rho -|z_n|)^2|w_n|^2},
\end{eqnarray*}
where $\Delta_\rho(z)$ is the disc of center $z$ and radius $\rho$. Hence, for any $r<R$, by taking $\rho=\frac{1}{2}(r+R)$ we infer
\begin{equation}
\Vert F_j\Vert _{ L^\infty (\Delta_{r}^n) }\leq \frac {2^nC_n^{\frac 1 2}R\epsilon_0}{\pi^{\frac n 2} (R-r)^n |w_n| }.
\end{equation}
Since $F_j(z',w_n)-f(z',w_n)=0$, $\forall z'\in\Delta_r^{n-1}$, we can write $F_j(z)=f(z)+(z_n-w_n)g_j(z)$ for some function $g_j(z)=\sum_{\alpha\in{\mathbb N}^n} a_{j,\alpha}z^{\alpha}$ on $\Delta_r^{n-1}\times\Delta_R$. By (1), we get
\begin{eqnarray*}
\Vert g_j\Vert _{ \Delta_{r}^n } = \Vert g_j\Vert _{ \Delta_{r}^{n-1}\times\partial\Delta_{r} }
&\leq& \frac {1} {r-|w_n|} \Big(\Vert F_j\Vert _{ L^\infty (\Delta_{r}^n) }+\Vert f\Vert _{ L^\infty (\Delta_{r}^n) }\Big)\\
&\leq& \frac {1} {r-|w_n|} \Big(\frac {2^nC_n^{\frac 1 2}R \epsilon_0} { \pi^{\frac n 2} (R-r)^n |w_n| }+\Vert f\Vert _{ L^\infty (\Delta_{r}^n) }\Big).
\end{eqnarray*}
Thanks to the Cauchy integral formula, we find
$$| a_{j,\alpha} |\leq \frac { \Vert g_j\Vert _{ \Delta_{r}^n } } { r^{ |\alpha | } }\leq \frac {1} { (r-|w_n|) r^{ |\alpha| } } \Big(\frac {2^nC_n^{\frac 1 2}R \epsilon_0} { \pi^{\frac n 2} (R-r)^n |w_n| }+\Vert f\Vert _{ L^\infty (\Delta_{r}^n) }\Big).$$
We take in any case $\eta\leq\epsilon_0\leq\epsilon\leq\frac{1}{2}r$. As $|w_n|<\eta\leq\frac{1}{2}r$, this implies
$$| w_n | | a_{j,\alpha} |\,r^{|\alpha |}\leq \frac {2}{r} \Big(\frac {2^nC_n^{\frac 1 2}R \epsilon_0} { \pi^{\frac n 2} (R-r)^n }+\Vert f\Vert _{ L^\infty (\Delta_{r}^n) } |w_n| \Big)\leq C'\epsilon_0,$$
for some constant $C'$ depending only on $n,\,r,\,R$ and~$f$. This yields the estimates of Lemma 2.1 for\break $\epsilon_0:=C''\epsilon$ with $C''$ sufficiently small.
Finally, we prove that $\IM(F_j)\leq \IM(f)$. Indeed, if $\IM(g_j)\geq\IM(f)$, since $|w_n\Vert a_{j,\alpha}|\leq r^{-|\alpha| }\epsilon$, we can choose $\epsilon$ small enough such that $\IM(F_j) = \IM(f)$ and $\Big|\frac {\textstyle\IC(F_j)}{\textstyle\IC(f) }\Big|\in (\frac 1 2, 2)$. Otherwise, if $\IM(g_j)<\IM(f)$, we have $\IM(F_j)=\IM(g_j)<\IM(f)$.\hfill$\square$
\vskip5pt
\n
{\it Proof of the main theorem.} By well-known properties of (pluri)potential theory, the $L^1$ convergence of $\psi$ to $\varphi$ implies that $\psi\to\varphi$ almost everywhere, and the assumptions guarantee that $\varphi$ and $\psi$ are uniformly bounded on every relatively compact subset of~$\Omega'$. In particular, after shrinking $\Omega'$ and subtracting constants, we can assume that $\varphi\le 0$ on~$\Omega$. Also, since the $L^1$ topology is metrizable, it is enough to work with a sequence $(\psi_j)_{j\geq 1}$ converging to $\varphi$ in $L^1(\Omega')$. Again, we can assume that $\psi_j\leq 0$ and that $\psi_j\to\varphi$ almost everywhere on $\Omega'$. By a trivial compactness argument, it is enough to show (i) and (ii) for some neighborhood $\Omega''$ of a given point $z_0\in\Omega'$. We assume here $z_0=0$ for simplicity of notation, and fix a polydisc $\Delta_R^n$ of center $0$ with $R$ so small that $\Delta_R^n\subset\Omega'$. Then \hbox{$\psi_j(\bu,z_n)\to\varphi (\bu,z_n)$} in the topology of $L^1(\Delta_R^{n-1})$ for almost every $z_n\in\Delta_R$.
\vskip3pt
\n
{\it Proof of statement {\rm(i)}}. We have here $\intn_{\Delta_R^n}e^{-2c\,\varphi}dV_{2n}<+\infty$ for $R>0$ small enough. By Lemma~2.1 with $f=1$, for every $r<R$ and $\epsilon>0$, there exist $w_n\in\Delta_{\epsilon}\smallsetminus\{0\}$, an index $j_0$, a number $\tilde c>c$ and a sequence of holomorphic functions $F_j$ on $\Delta_r^n$, $j\geq j_0$, such that $F_j(z)=1+(z_n-w_n)\sum a_{j,\alpha}z^\alpha$, $|w_n||a_{j,\alpha}|\,r^{-|\alpha|}\le\epsilon$ and
$$\intn_{\Delta_r^n}|F_j(z)|^2e^{-2\tilde c\,\psi_j(z)}dV_{2n}(z)\le\frac{\epsilon^2}{|w_n|^2},~~~\forall j\geq j_0.$$
For $\epsilon\leq\frac{1}{2}$, we conclude that $|F_j(0)|=|1-w_na_{j,0}|\geq\frac{1}{2}$ hence $c_{\psi_j}(0)\geq\tilde c>c$ and the first part of (i) is proved. In fact, after fixing such $\epsilon$ and $w_n$, we even obtain the existence of a neighborhood $\Omega''$ of $0$ on which $|F_j|\geq\frac{1}{4}$, and thus get a uniform bound $\int_{\Omega''}e^{-2\tilde c\,\psi_j(z)}dV_{2n}(z)\le M<+\infty$. The second assertion of (i) then follows from the estimate
\begin{eqnarray*}
\intn_{\Omega''}\big|e^{-2c\,\psi_j}-e^{-2c\,\varphi}\big|
dV_{2n}&\leq&
\intn_{\Omega''\cap\{|\psi_j|\leq A\}}\big|e^{-2c\,\psi_j}-e^{-2c\,\varphi}\big|dV_{2n}\\
&+&\intn_{\Omega''\cap\{\psi_j<-A\}}e^{-2c\,\varphi}dV_{2n}\\
&+&e^{-2(\tilde c-c)A}\intn_{\Omega''\cap\{\psi_j<-A\}}e^{-2\tilde c\,\psi_j}dV_{2n}.
\end{eqnarray*}
In fact the last two terms converge to $0$ as $A\to+\infty$, and, for $A$ fixed, the first one in the right hand side converges to $0$ by Lebesgue's bounded
convergence theorem, as $\psi_j\to\varphi$ almost everywhere on
$\Omega''$.
\vskip3pt

\n{\it Proof of statement {\rm (ii)}.} Take $f_1,\ldots,f_k\in\mathcal O_{\mathbb C^n, 0}$ such that $(f_1,\ldots,f_k)$ is a standard basis of ${\mathcal I}(c\,\varphi)_0$ with $\IM(f_1)<\ldots<\IM(f_k)$, and $\Delta_R^n$ a polydisc so small that
$$\intn_{\Delta_R^n}|f_l(z)|^2e^{-2c\,\varphi(z)}dV_{2n}(z)<+\infty,~~~l=1,\ldots,k.$$
Since the germ of $f$ at $0$ belongs to the ideal $(f_1,\ldots,f_k)$, we can essentially argue with the $f_l$'s instead of~$f$. By Lemma~2.1, for every $r<R$ and $\epsilon_l>0$, there exist $w_{n,l}\in\Delta_{\epsilon_l}\smallsetminus\{0\}$, an index $j_0=j_0(w_{n,l})$, a number $\tilde c=\tilde c(w_{n,l})>c$ and a sequence of holomorphic functions $F_{j,l}$ on $\Delta_r^n$, $j\geq j_0$, such that $F_{j,l}(z)=1+(z_n-w_{n,l})\sum a_{j,l,\alpha}z^\alpha$, $|w_{n,l}||a_{j,l,\alpha}|\,r^{-|\alpha|}\le\epsilon_l$ and
\begin{equation}
\intn_{\Delta_r^n}|F_{j,l}(z)|^2e^{-2\tilde c\,\psi_j(z)}dV_{2n}(z)\le\frac{\epsilon_l^2}{|w_{n,l}|^2},~~~\forall l=1,\ldots,k,~~\forall j\geq j_0.
\end{equation}
Since $\psi_j\leq\varphi$ and $\tilde c>c$, we get $F_{j,l}\in{\mathcal I}(\tilde c\,\psi_j)_0\subset {\mathcal I}(c\,\varphi)_0$. The next step of the proof consists in modifying $(F_{j,l})_{1\leq l\leq k}$ in order to obtain a standard basis of ${\mathcal I}(c\,\varphi)_0$. For this, we proceed by selecting successively $\epsilon_1\gg\epsilon_2\gg\ldots\gg\epsilon_k$ (and suitable $w_{n,l}\in\Delta_{\epsilon_l}\smallsetminus\{0\}$). We have $\IM(F_{j,1}),\ldots,\IM(F_{j,k})\in\IM(\mathcal I(c\,\varphi)_0)$, in particular $\IM(F_{j,1})$ is divisible by $\IM(f_l)$ for some $l=1,\ldots,k$. Since $\IM(F_{j,1})\leq\IM(f_1)<\ldots<\IM(f_k)$, we must have $\IM(F_{j,1})=\IM(f_1)$ and thus $\IM(g_{j,1})\geq\IM(f_1)$. As $|w_{n,1}||a_{j,1,\alpha}|\le\epsilon_1$, we will have $\Big|\frac {\textstyle\IC(F_{j,1}) } {\textstyle\IC(f_1) }\Big|\in (\frac 1 2, 2)$ for $\epsilon_1$ small enough. Now, possibly after changing $\epsilon_2$ to a smaller value, we show that there exists a polynomial $P_{j,2,1}$ such that the degree and coefficients of $P_{j,2,1}$ are uniformly bounded, with $\IM(F_{j,2}-P_{j,2,1} F_{j,1}) = \IM(f_2)$ and $\frac {\textstyle|\IC (F_{j,2}-P_{j,2,1} F_{j,1})|} {\textstyle|\IC(f_2)|}\in (\frac 1 2, 2)$. We consider two cases:
\vskip5pt
\n
{\bf Case 1}\/: If $\IM(g_{j,2})\geq\IM(f_2)$, since $|w_{n,2}||a_{j,2,\alpha}|\leq r^{-|\alpha|}\epsilon_2$, we can choose $\epsilon_2$ so small that $\IM(F_{j,2}) = \IM(f_2)$ and $\frac {\textstyle|\IC(F_{j,2})|}{\textstyle|\IC(f_2)|}\in (\frac 1 2, 2)$. We then take $P_{j,2,1}=0$.
\vskip5pt
\n
{\bf Case 2}\/: If $\IM(g_{j,2})<\IM(f_2)$, we have $\IM(g_{j,2}) = \IM(F_{j,2})\in \IM(\mathcal I(c\,\varphi)_0)$. Hence $\IM(g_{j,2})$ is divisible by $\IM(f_l)$ for some $l=1,\ldots,k$. However, since $\IM(g_{j,2})<\IM(f_2)<\ldots<\IM(f_k)$, the only possibility is that $\IM(g_{j,2})$ be divisible by $\IM(f_1)$. Take $b\in\mathbb C$ and $\beta,\gamma\in{\mathbb N}^n$ such that $\IT(g_{j,2}):=a_{j,2,\gamma}z^{\gamma}=bz^{\beta}\,\IT(F_{j,1})$. We have $z^{\beta}\leq z^{\gamma}=\IM(g_{j,2})<\IM(f_2)$ and
$$|w_{n,2}||b|=|w_{n,2}|\frac {|\IC(g_{j,2})|} {|\IC(F_{j,1})|}\leq
\frac {2|w_{n,2}||a_{j,2,\gamma}|} {|\IC(f_1)|}
\leq \frac{2r^{-|\gamma|}\epsilon_2}{|\IC(f_1)|}$$
can be taken arbitrarily small.
Set $\tilde g_{j,2}(z)=g_{j,2}(z)- bz^{\beta} F_{j,1}(z)=\sum \tilde a_{j,2,\alpha} z^{\alpha}$ and
$$\tilde F_{j,2}(z)=f_2(z)+(z_n-w_{n,2})\tilde g_{j,2}(z)=F_{j,2}(z)-b(z_n-w_{n,2})z^{\beta} F_{j,1}(z).$$
We have $\IM(\tilde g_{j,2}) > \IM(g_{j,2})$. Since $|w_{n,2}||b|=O(\epsilon_2)$ and $|w_{n,2}|| a_{j,2,\alpha}|=O(\epsilon_2)$, we get $|w_{n,2}||\tilde a_{j,2,\alpha}|=O(\epsilon_2)$ as well. Now, we consider two further cases. If $\IM(\tilde g_{j,2}) \geq \IM(f_2)$, we can again change $\epsilon_2$ for a smaller value so that $\IM(\tilde F_{j,2}) = \IM(f_2)$ and $\frac {\textstyle|\IC (\tilde F_{j,2})|}{\textstyle|\IC (f_2)|}\in (\frac 1 2, 2)$. Otherwise, if $\IM(\tilde g_{j,2}) < \IM(f_2)$, we have $\IM(F_{j,2}) = \IM(g_{j,2}) < \IM( \tilde F_{j,2} ) = \IM( \tilde g_{j,2} ) < \IM(f_2)$. Notice that \hbox{$\{ z^{ \gamma }:\ z^{ \gamma }<\IM(f_2) \}$} is a finite set. By using similar arguments a finite number of times, we find $\epsilon_2$ so small\break that $\IM(F_{j,2}-P_{j,2,1}F_{j,1}) = \IM(f_2)$ and $\frac {\textstyle|\IC(F_{j,2}-P_{j,2,1} F_{j,1})|}{\textstyle|\IC(f_2)|}\in (\frac 1 2, 2)$ for some polynomial $P_{j,2,1}$.
\n
Repeating the same arguments for $F_{j,3},\ldots,F_{j,k}$, we select inductively $\epsilon_l$, $l=1,\ldots,k$, and construct linear combinations $F'_{j,l}=F_{j,l}-\sum_{1\leq m\leq l-1}P_{j,l,m}F'_{j,m}$ with polynomials $P_{j,l,m}$, $1\leq m<l\le k$, possessing uniformly bounded coefficients and degrees, such that $\IM(F'_{j,l})=\IM(f_l)$ and
$\frac {\textstyle|\IC (F'_{j,l})|} {\textstyle|\IC(f_l)|}\in (\frac 1 2, 2)$ for all $l=1,\ldots,k$ and $j\ge j_0$. This implies that
$(F'_{j,1},\ldots,F'_{j,k})$
is also a standard basis of $\mathcal I(c\,\varphi)_0$. By Theorem 1.2.2 in \cite{Ga79}, we can find $\rho$,~\hbox{$K>0$} so small that there exist holomorphic functions $h_{j,1},\ldots,h_{j,k}$ on $\Delta_{\rho}^n$ with $\rho<r$, such that
$$f=h_{j,1}F'_{j,1}+h_{j,2}F'_{j,2}+\ldots+h_{j,k}F'_{j,k}~~~
\hbox{on $\Delta_{\rho}^n$}$$
and $\Vert h_{j,l}\Vert _{L^\infty (\Delta_{\rho}^n)}\leq K\Vert f\Vert _{L^\infty (\Delta_{r}^n)}$, for all $l=1,\ldots,k$ ($\rho$ and $K$ only depend on $f_1,\ldots,f_k$). By (2), this implies a uniform bound
$$\intn_{ \Delta_\rho^n } |f(z)|^2e^{-2\tilde c\,\psi_j(z)} d V_{2n}(z)\le M<+\infty$$
for some $\tilde c>c$ and all $j\geq j_0$. Take $\Omega''=\Delta_\rho^n$.
We obtain the $L^1$ convergence of $|f|^2e^{-2c\,\psi_j}$ to $|f|^2e^{-2c\,\varphi}$ almost exactly as we argued for the second assertion of part~(i), by using the estimate
\begin{eqnarray*}
\intn_{\Omega''}|f|^2\big|e^{-2c\,\psi_j}-e^{-2c\,\varphi}\big|
dV_{2n}&\leq&
\intn_{\Omega''\cap\{|\psi_j|\leq A\}}|f|^2\big|e^{-2c\,\psi_j}-e^{-2c\,\varphi}\big|dV_{2n}\\
&+&\intn_{\Omega''\cap\{\psi_j<-A\}}|f|^2e^{-2c\,\varphi}dV_{2n}\\
&+&e^{-2(\tilde c-c)A}\intn_{\Omega''\cap\{\psi_j<-A\}}|f|^2e^{-2\tilde c\,\psi_j}dV_{2n}.
\end{eqnarray*}
\vskip3pt
\n
{\bf Acknowledgments:} The author is deeply grateful to Professor Jean-Pierre Demailly for valuable comments and helpful discussions during the preparation of this work. The research was done while the author was supported by the ANR project MNGNK, decision N$^\circ$ ANR-10-BLAN-0118. He would like to thank Professor Andrei Teleman and the members of the LATP, CMI, Universit\'e de Provence, Marseille, France for their kind hospitality.

\n
Pham Hoang Hiep

\n
Department of Mathematics, Hanoi National University of Education

\n
136-Xuan Thuy, Cau Giay, Hanoi, Vietnam

\n
{\it e-mail\/}: {\tt phhiep$_-$vn@yahoo.com}

\end{document}